\documentclass{article}
\usepackage{etex}
\usepackage[usenames,dvipsnames,svgnames,table]{xcolor}
\usepackage{color}
\usepackage[latin1]{inputenc}
\usepackage{amsmath, enumerate, amsfonts}
\usepackage[widemargins]{a4}
\usepackage{tikz} 
\usepackage{tikz-cd}  
\usepackage{amsmath, enumerate, amsfonts,amssymb, pictexwd, amscd}

\newcommand{\pcom}{^\wedge_p}
\usepackage[all]{xy}
\xyoption{arc}
\newtheorem{Conjecture}{Conjecture}[section]
\newtheorem{Corollary}[Conjecture]{Corollary}
\newtheorem{Theorem}[Conjecture]{Theorem}

\newtheorem{Proposition}[Conjecture]{Proposition}
\newtheorem{Lemma}[Conjecture]{Lemma}
\newtheorem{Remark}[Conjecture]{Remark}

\usepackage[latin1]{inputenc}
\usepackage{graphicx}
\usepackage{amsmath, enumerate, amsfonts}
\usepackage[widemargins]{a4}
\usepackage[all]{xy}
\xyoption{arc}
\color{black}

\title{A few examples of $p-$good and $p-$bad classifying spaces}
\author{Nora Seeliger}
\date{}
\begin{document}
\maketitle
\begin{abstract}
We give examples of spaces which are good and bad at different primes in the sense of Bousfield and Kan \cite{BK} in any arbitrary combination and investigate which impact the existence of a Sylow $p$-subgroup has on the homotopy type on the classifying space and under which conditions the homotopy type of wedges of classifying spaces is good or bad for a solid ring $R$. We give results relating to various other $R$-homological structures and a collection of examples.
\end{abstract}
\section{Introduction}
In general, it is a very hard question to decide whether a space is good or bad in the sense of Bousfield and Kan [\textit{Homotopy Limits, Completions and Localizations}, Springer Lecture Notes in Mathematics, Springer-Verlag, Berlin, Heidelberg, New York, 1972] and there are not many explicit examples in the literature. In this note we give results about $R-$completed classifying spaces for various solid rings $R$. In particular, we give examples of spaces which are good at some primes and bad at other primes in any arbitrary combination. The special case answers a question raised during Bob Oliver's talk with the title "Local structure of groups and of their classifying spaces" at the MSRI Berkeley at the beginning of the algebraic topology program on January 30th 2014: Are there examples of spaces which are good for some primes and bad for other primes simultaneously? This is answered affirmatively for arbitrary combinations of sets of primes by providing a collection of explicit examples of such spaces. 
 The general statement that a space is good for one solid ring and bad for a different solid ring is false as shown in a separate manuscript \cite{RgoodRbad}. We continue to illustrate the impact of the existence of a finite Sylow $p$-subgroup on the homotopy type of the classifying space of a group. 
 Our results on completion of spaces apply to proving spaces not to be retracts of others in the case of certain projective spaces. We give an argument that illustrates the limitation of the mod-$R$ fibre Lemma to show that a space is bad for subrings of the rationals and some results relating the $R-$completion and $R$-localisation of a space to various other $R$-homological structures including the structure of its $\mathbb{F}_p-$
homology structure as an unstable coalgebra over the dual Steenrod algebra and a collection of examples. There are also applications to the homotopy type of group models for fusion systems \cite{Robinson1}, \cite{Ian+Radu} which were first studied by ourselves in \cite{gmffs}.\\
\label{cIntro} The author was supported by ARC Discovery Project DP120101399 at the Australian National University, the Mathematical Sciences Research Institute in Berkeley, a postdoctoral fellowship at the University of Haifa, and at the Friedrich-Schiller-Universit\"at Jena. I would like to thank David Benson, David Blanc, Bill Dwyer, Emmanuel Dror Farjoun, Achim Krause, Assaf Libman, and Bob Oliver for enriching discussion on the subject and the Mathematical Sciences Research Institute in Berkeley for its hospitality.
\section{Preliminaries}
The $p$-completion in the sense of Bousfield-Kan \cite{BK} is a special case of the $R-$completion functor of Bousfield and Kan which is defined for solid rings $R$. This is a functor from the category of simplicial sets to itself. The category of simplicial sets will be denoted by spaces thoughout this article. A ring $R$ is \textbf{solid} if the map $R\otimes _{\mathbb{Z}}R\rightarrow R$ is an isomorphism. Solid rings are $\mathbb{Z}/n\mathbb{Z}$, subrings of the rationals $\mathbb{Z}[J^{-1}]$ for any set $J$ of primes, the product rings $\mathbb{Z}[J^{-1}]\times \mathbb{Z}/n\mathbb{Z}$ where each prime factor of $n$ is in $J$, and direct limits of these three types of rings \cite{BK2}. The Bousfield-Kan completion is related to completions and localizations in the sense of Malcev \cite{Malcev}, Sullivan \cite{Sullivan}, Quillen \cite{Quillen}.
 In $p-$local homotopy theory the
classifying space  $|\mathcal{L}|\pcom$ of a $p-$local finite group $(S,\mathcal{F},\mathcal{L})$ is one of the main 
objects of study where $(-)\pcom$
denotes the $p-$completion functor with respect to the ring $\mathbb{Z}/p\mathbb{Z}$.  A space X is called \textbf{$p-$complete} if the map $X\rightarrow X\pcom$ is a mod $p-$equivalence and  \textbf{$p-$good} if the natural map $H_*(X;\mathbb{F}_p)\rightarrow H_*(X\pcom ;\mathbb{F}_p)$ is an isomorphism. Otherwise the space is called \textbf{$p$-bad}. The properties good and bad are defined analogously for all solid rings. Examples of spaces which are good for all primes and the integers are classifying spaces of finite groups. Spaces with finite homotopy groups are good for all solid rings $R$. For spaces with finite homotopy groups or finite homology groups the $\mathbb{Z}-$completion is up to homotopy, the product of the $\mathbb{Z}/p\mathbb{Z}-$completions. The space $\mathbb{R}P^2$ is good for all primes and bad for $\mathbb{Z}_{(J)}$ as long as $2\in J$, and bad for the integers \cite{BK}. In a different manuscript we show that the space $S^1\vee S^1$ is bad for all nonzero solid rings \cite{RgoodRbad}.
 A group $\pi$ acts on a group $G$
if there is given a homomorphism
$\alpha :\pi\rightarrow Aut(G)$
and such an action is called \textbf{nilpotent} if there exists a finite
sequence of subgroups of $G$:
$G=G_1\supset\cdots\supset G_j\supset\cdots\supset G_n=1$
such that for each $j$ we have
$G_j$ is closed under the action of $\pi$,
$G_{j +1}$ is normal in $G_j$, and $G_j/G_{j +1}$ is abelian, and
the induced action on $ G_j/G_{j+1} $ is trivial.
A (possibly infinite) group $G$ has a \textbf{Sylow} $p-$subgroup $S$ if it has a subgroup isomorphic to $S$ and all other subgroups of $G$ whose order is a power of $p$ are subconjugate to $S$. Examples of infinite groups with finite Sylow $p-$subgroups are group models for fusion systems in the sense of Robinson and Leary-Stancu and ourselves and ourselves together with Libman. Subgroups of amalgamated products are described by Kurosh's subgroup Theorem. The most general version states that \cite{Massey} if $H$ is a subgroup of the free amalgamated product of groups $ \underset{i\in I}{*}G_i = G$, then $H=F(X)*(*_{j\in J} g_jH_jg_j^{-1})$, where $X$ is a subset  of $G$ and $J$ is some index set and $g_j \in G$ and each $H_j$ is a subgroup of some $G_i$. Recall that for any pair of groups $ G,H$  we have a weak equivalence of classifying spaces $B(G*H)\simeq BG\vee BH$, \cite{Massey}. The mod $R-$Fibre Lemma \cite[Lemma 5.1]{BK} states that the $R$-completion preserves fibrations up to homotopy of connected spaces $F\rightarrow E\rightarrow B$ for which the fundamental group of the base space $\pi _1(B)$ acts nilpotently on every reduced homology group of the fibre with coefficients in the ring of the completion $\overline{H}_i(F;R)$. It is used to show that spaces are good or bad without constructing the completion tower. 
 A space is \textbf{nilpotent} if the fundamental group acts nilpotently on all higher homotopy groups \cite{BK}. A group $G$ is $R-$\textbf{nilpotent} if it has a finite central series $G=G_1\subset\cdots G_j\cdots G_n=1 $ such that each quotient $G_j/G_{j+1}$ admit an $R$-module structure. A map $f:\{G_s\}\rightarrow \{H_s\}$ between two towers of groups is a \textbf{pro-isomorphism} if, for every group $B$, it induces an isomorphism 
\begin{eqnarray}\underset{\leftarrow}{lim}\text{ }Hom_{(groups)}(H_s,B)\cong \underset{\leftarrow }{lim}\text{ }Hom_{(groups)}(G_s,B).\end{eqnarray}
The completion tower $R_nX$defined by Bousfield and Kan preserves the homology with $R$-coefficients \cite{Farjoun2}. This property of preserving $R$-homology characterizes the tower completely \cite{BK}. 
The general statement whether a space can be $R$-good and $R'$-bad for arbitrary solid rings $R\neq R'$ is false \cite{RgoodRbad}.
\section{On nilpotent group actions and the mod $R-$Fibre Lemma}
We start by giving some results concerning nilpotent group actions and the mod $R-$Fibre Lemma.
\begin{Remark}
For a prime $p$ a
 finite $p-$group always acts nilpotently on a finite $p-$group \cite[p.47]{Hall}.
\end{Remark}
We illustrate that this result does not generalize to free amalgamated products of finite $p-$groups.
\begin{Lemma}
 The group $C_p*C_p$ does not act nilpotently on the finite $p-$group $C_p\times C_p$ if each factor of the direct product is embedded as upper and lower triangular matrices respectivly in $SL(2,p)$.
\end{Lemma}
\underline{Proof:} 
 Assume not. The group $SL(2,p)$ acts transitively on lines in $\mathbb{F}_p^2$ and is generated by any two Sylow $p-$subgroups. So it acts nontrivially on $C_p\times C_p$ and the various embeddings of $C_p$ \cite{Dickson}. 
$\Box$
\begin{Corollary}
Let $G$ be a symmetric group which acts nonnilpotently on a finite $2$-group $H$. Then there is a free amalgamated product of copies of $C_2$ which acts nonnilpotently on $H$.
\end{Corollary}
\underline{Proof:} The symmetric group $G$ is generated by elements of order $2$. $\Box$

The following proposition illustrates that we can apply the mod$-R$ fibre lemma \cite[Lemma 5.1]{BK} in $\mathbb{F}_p-$homology to any fibration where the fundamental group of the base space is a finite $p-$group.
\begin{Proposition}
A finite $p-$group $P$ acts nilpotently on an $\mathbb{F}_p[P]-$module of arbitrary cardinality.
\end{Proposition}
\label{nilpotent}
\underline{Proof:} Let $I<\mathbb{F}_p[P]$ be the augmentation ideal and denote the module by $M$. If $n>0$ is such that $I^n\ne0$, then $I^{n+1}<I^n$. This implies that there is some $k$ such that $I^k=0$. Since $P$ acts trivially on each quotient $I^nM/I^{n+1}M$ of the sequence of submodules $M>IM>I^2M>....>I^kM=0$, the module $M$ is nilpotent.$\Box$\\[0.3cm]

\begin{Remark}
Alternatively the statement of the preceeding proposition can be proved as follows.
\end{Remark}
\begin{Proposition}
Let $G$ be a virtually finite group and denote $G\pcom$ the group completion of $G$, by which we mean the inverse limit of the respective quotients of the group over its $p$-lower central series as defined by Bousfield, see \cite{Bousfield}. Then we have a weak equivalence of spaces $(BG)\pcom\simeq B(G\pcom)$.
\end{Proposition}
\underline{Proof:} This follows from the construction of the completion functor \cite[Chapter 4]{BK}. $\Box$

\begin{Proposition}
Let $G$ be a virtually finite $p$-group. Then the classifying space BG is p-good iff and only if the index of the free subgroup is trivial or one.
\end{Proposition}
\underline{Proof:} This follows analogously to the argument in the previous proof. $\Box$\\[0.3cm]
\begin{Proposition}
Let $p$ be a prime and $\{G_i\}_{i\in I}$ a family of finite groups with $(|G_i|,p)=1$ for all $i\in I$. Let $G$ be an arbitrary combination of direct and amalgamated product of the $G_i$. Assume that the group $G$ has a finite index normal subgroup generated by elements of $q'-$order. Then classifying space of the group $G$ is $p$-good.  
\end{Proposition}
\underline{Proof:} The group $G$ has a finite index normal subgroup generated by elements of $q'-$order and therefore \cite[Proposition 1.11]{AKO} the classifying space is $q-$good. The proof of the result for the classifying space of an arbitrary combination of direct and amalgamated products of the groups $G_{p_i}$ is totally analogous. $\Box$\\[0.3cm]
The following lemma allows to construct bad spaces for arbitrary solid rings $R$ and will be used at a later point. 
 \begin{Lemma}
 \label{pbad}
 Let $A$ and $B$ be spaces such that $A$ is $R$-bad for a solid ring $R$ and $A$ is a retract of $B$. Then the space $B$ is $R$-bad.
 \end{Lemma}
 \underline{Proof:} This follows directly from the definition of $R$-badness and it being a universal property.$\Box$
The mod $R-$Fibre lemma and our result on nilpotent group actions Proposition \ref{nilpotent} allow to prove.
\begin{Theorem}
Let $X$ be a space which fibers over a connected space with fundamental group a finite $p$-group. Then $X$ is $p-$good iff the fibre $F$ is $p-$good.
\end{Theorem}
\underline{Proof:} The fundamental group of the base acts nilpotently on the reduced homology of the fibre and therefore we can $p$-complete the fibration. So the fibre is $p$-good if and only if the total space is $p$-good.$\Box$
\begin{Proposition}
\label{pgoodpbad}
Let $p$ be a prime. Then $ B(C_p*C_p)\simeq BC_p\vee BC_p$ is $p$-good if and only if $p=2$.
\end{Proposition}
\underline{Proof:} There is a canonical map $C_p*C_p\rightarrow C_p\times C_p$ for all primes $p$. It
follows from Kurosh's subgroup theorem that the kernel is a free group $\Gamma$
and an Euler characteristic computation in the sense of C.\ T.\ C.\ Wall which is referenced in Serre \cite{Serre} gives that $\Gamma $ has rank $(p-1)^2$. 
So for
$p$ odd 
this rank is strictly bigger than one. A finite $p$-group always acts nilpotently on a finite $p-$group. 
Therefore $C_p\times C_p$ acts nilpotently on the $\mathbb{F}_p$-homology of $\Gamma$ and therefore we can $p$-complete the fibration $B\Gamma\rightarrow BC_p\vee BC_p\rightarrow BC_p\times BC_p$. Therefore $BC_p\vee BC_p$ is $p$-bad since if it was $p$-good the space $B\Gamma$ would be $p$-complete which means it would be $p$-good which is a contradiction. 
For $p=2$ we have that $C_2*C_2$ is isomorphic to the infinite dihedral group $D_{\infty}$. It follows from the mod-$R$ fibre Lemma \cite[Lemma 5.1]{BK} for the fibration $B\Gamma\rightarrow BD_{\infty}\rightarrow C_2\times C_2$ that $BD_{\infty}$ is $2$-good
because the rank of $\Gamma $ is $1$ and the space $B\Gamma\simeq
S^1$ is good for all primes, so in particular
good for the prime $2$. Therefore
 $B(C_2*C_2)\simeq BC_2\vee BC_2$ is $2$-good.
 $\Box$
\begin{Remark}
In this case the completions of Bousfield Kan \cite{BK}, Malcev \cite{Malcev}, Quillen \cite{Quillen}, Sullivan \cite{Sullivan} coincide.
\end{Remark}
More generally we have.
\begin{Remark}
The spaces $\mathbb{R}P^k$ are $2$-good for all $k\geq 0$. This includes the possibility that $k=\infty$.
\end{Remark}
\underline{Proof:} We have a fibration up to homotopy $S^k\rightarrow \mathbb{R}P^k\rightarrow K(C_2,1)$. The action of the fundamental group on the $\mathbb{F}_2$-homology of the fibre is nilpotent. Therefore we can $2$-complete the fibration. Since the spaces $S^k$ and $K(C_2,1)$ are both $2-$good the space $\mathbb{R}P^k$ is $2$-good. Recall that $\mathbb{R}P^{\infty}\simeq BC_2$.$\Box$
 \begin{Proposition}
\label{groups}
Let $\{P_i\}_{i\in I}$ be an arbitrary collection of $p-$groups.
The classifying space $B(\underset{i\in I}{*}P_i)$ of the free amalgamated product is $p-$bad as long as the set $I$ has at least two elements for $p$ odd and at least three elements for $p=2$ or not both groups are isomorphic to the cyclic group $C_2$. 
 \end{Proposition}
\underline{Proof:} This follows analogously to the proof of Proposition \ref{pgoodpbad} and together with Lemma \ref{pbad}. $\Box$
\begin{Proposition}
\label{pgoodpbad2}
Let $p$ be a prime. Then $ B(C_p*C_p)\simeq BC_p\vee BC_p$ is $p$-good if and only if $p=2$.
\end{Proposition}
\underline{Proof:} We can inspect the lower $p$-central series of the group $C_p*C_p$ for all primes $p$. Look at the canonical group homomorphism $C_p * C_p \rightarrow C_p \times C_p$. The kernel is a free group. An Euler characteristic computation shows that this group has rank $(p-1)^2$. The action of the quotient on the cohomology
of the kernel is nilpotent, so the mod $R-$Fibre Lemma applies to the induced map on classifying spaces. So when we complete, the kernel 
becomes a free $p$-nilpotent group on $(p-1)^2$ generators $K$. This
group has uncountable cohomology for $p>2$, and for placement reasons a lot of it survives in the Serre
spectral sequence for the group extension $K\rightarrow (C_p*C_p)\pcom\rightarrow (C_p\times C_p)\pcom$. 
\begin{Proposition}
Let $G_1,G_2$ be groups, not both isomophic to $C_2$, and $p$ a prime. In the case that the classifying spaces $BG_1$, $BG_2$ are both $p-$good assume that the natural action of the product $G_1\times G_2$ is nilpotent on the $\mathbb{F}_p-$homology of the fibre of the fibration $B(G_1* G_2)\rightarrow B(G_1\times G_2)$. Then the space $BG_1\vee BG_2$ is $p-$bad.
\end{Proposition}
\underline{Proof:} In the case where Lemma \ref{pbad} does not apply this follows directly from the mod $R-$fibre Lemma \cite[Lemma 5.1]{BK} and the associated homology spectral sequence. Note that after Kurosh's subgroup Theorem the fibre has the homotopy type of a finite wedge of $1-$spheres and the reduced $\mathbb{F}_p-$homology of the fibre is isomorphic as group to a finite number of copies of a cyclic group.$\Box$

\begin{Proposition}
Let $G_1,G_2$ be finite groups, not both isomophic to $C_2$, and $p$ a prime. In the case that the classifying spaces $BG_1$, $BG_2$ are both $p-$good assume that the natural action of the product $G_1\times G_2$ is nilpotent on the $\mathbb{F}_p-$homology of the fibre of the fibration $B(G_1* G_2)\rightarrow B(G_1/O^p(G_1)\times O^p(G_2))$ and that $G_1/O^p(G_1)$ and $G_1/O^p(G_1)$ are both nontrivial, where $O^p(G)$ denotes the normal subgroup of $G$ generated by elements whose order does not get divided by the prime $p$. Then the space $BG_1\vee BG_2$ is $p-$bad.
\end{Proposition}
\underline{Proof:} In the case where Lemma \ref{pbad} does not apply this follows directly from the mod $R-$fibre Lemma \cite[Lemma 5.1]{BK} and the associated homology spectral sequence. Note that after Kurosh's subgroup Theorem the fibre has the homotopy type of a finite wedge of $1-$spheres and the reduced $\mathbb{F}_p-$homology of the fibre is isomorphic as group to a finite number of copies of a cyclic group.$\Box$

 \begin{Proposition}
Let $p$ be a prime and $P$ a finite nontrivial
$p-$group.
Then the
classifying
space of the free amalgamated product of $P$ with the integers 
$B(\mathbb{Z}*P)\simeq S^1\vee BP$ is bad for the prime $p$.
\end{Proposition}
\underline{Proof:}
There is the canonical projection map $\mathbb{Z}*P
\rightarrow P$. Denote the kernel of this map $\Gamma$. It follows from Kurosh's subgroup theorem that $\Gamma$ contains as a freely amalgamated factor several freely amalgamated copies of the integers. Therefore the classifying space $B\Gamma$ is weakly equivalent to a wedge product of a space $X$ and a bouquet of at least two spheres. Therefore the classifying space $B\Gamma$ is $q
-$bad
for the
 prime $p$.
 The group $ P
$ acts nilpotently
 on the $\mathbb{F}_p$-homology of $B\Gamma$ and therefore we can $p$-complete the fibration $B\Gamma\rightarrow S^1
\vee BP
\rightarrow S^1
\times BP$. Therefore $S^1\vee BP$
 is $p$-bad since if it was $p$-good the space $B\Gamma$ would be $p$-complete which means it would be $p$-good which is a contradiction. $\Box$\\[0.3cm]

The following proposition is anologue of the K\"unneth theorem.
\begin{Lemma}
\label{pgood}
Let $X$ and $Y$ be spaces and $R$ a solid ring such that $R$ is a principal domain, $X$ and $Y$ are both $R$-good and $Tor(H_p(X;R),H_q(Y;R))=0$ for all positive integers $p,q$ such that $p+q\geq 1$. The space $X\times Y$ is $R-$good.
\end{Lemma}
\underline{Proof:} The K\"unneth theorem is an isomorphism in this case for the trivial $R-$module $R$ itself. $\Box$\\[0.2cm]
As a special case we obtain the following. 
\begin{Remark}
This means the classifying space for a Leary-Stancu model \cite{Ian+Radu}
 can be $p-$good for as long as it consists of a single and not an iterated $HNN-$extension.
\end{Remark}
\section{The impact of Sylow $p$-subgroups on the homotopy type}
All group models involving only amalgams of finite groups have classifying spaces which are $p$-good.\\[0.3cm] 
The existence
of a
Sylow $p$-subgroup has implications on
the homotopy type of 
a
classifying space.
\begin{Remark}
 The classifying space of the group  $C_p*C_p$ is $p-$bad for $p$ odd as proved in Proposition \ref{pgoodpbad}.
\end{Remark}
 
\begin{Proposition}
The classifying space $BG$ of an iterated HNN construction $G=\frac{S*F[t_1,t_2,\cdots]}{<t_iut_i^{-1}=\phi _i(u)>}$ is bad for all nonzero solid rings regardless of whether $t_1, t_2, ...$ is generated by finitely many elements or infinitely. 
\end{Proposition}
\underline{Proof:} 
The classifying space $BG$ contains
 the space $S^1\vee S^1$ as a retract as the classifying space of the double iterated HNN construction over the trivial group.
 The space $S^1\vee S^1$
 is bad for all nonzero solid rings as proved by us in a separate manuscript \cite{RgoodRbad}.
$\Box$
 \begin{Proposition}
 Let $p$ be a prime and let $G_1, G_2$ be groups generated by elements of finite order and with finite Sylow $p-$subgroups $S_1$, $S_2$ respectively such that $p\neq2$ or not both $G_1$ and $G_2$ are isomorphic to $C_2$. Then $BG_1\vee BG_2$ is $p$-good if and only if either $S_1$ or $S_2$ or both are trivial.
 \end{Proposition}
\underline{Proof:} As above we have $BG_1\vee BG_2\simeq B(G_1*G_2)$.  If $p$ does not divide one of $|S_1|$ and $|S_2|$ then 
\cite[ Proposition 1.11]{AKO}
 applies to show that 
$BG_1\vee BG_2$ is $p-$good. In the remaining case that $p$ divides both $|G_1|$ and $|G_2|$ it follows from Kurosh's subgroup theorem that the canonical map $G_1*G_2\rightarrow G_1\times G_2/<g_i| g_i\in G_i, i\in\{1,2\},g_i \text{ has } p'-\text{order}>$ has kernel $\Gamma$ which is isomorphic to a free product on more than one generator and free groups generated by elements of finite and $p'-$order. The quotient $ G_1\times G_2/<g_i| g_i\in G_i, i\in\{1,2\},g_i \text{ has } p'-\text{order}>$ is  a finite $p-$group and therefore acts nilpotently on the $\mathbb{F}_p$-homology of $\Gamma$ and therefore we can $p$-complete the fibration $B\Gamma\rightarrow BG_1\vee BG_2\rightarrow B( G_1\times G_2/<g_i| g_i\in G_i, i\in\{1,2\},g_i \text{ has } p'-\text{order}>)$. Therefore $BG_1\vee BG_2$ is $p$-bad since if it was $p$-good the space $B\Gamma\pcom$ would be $p$-complete which means it would be $p$-good which is a contradiction. $\Box$
\begin{Proposition}
Let $G$ be a group with a finite Sylow $p-$subgroup $S$. Let $H$ be the subgroup of $G$ generated by all elements whose order does not get divided by $p$ and by all elements of infinite order. Then the classifying space $BG$ is $p-$good if and only if the classifying space $BH$ is $p-$good.
\end{Proposition}
\underline{Proof:} Note that $H\unlhd  G$ and consider the fibration $BH\rightarrow BG\rightarrow B(G/H)$. Applying the mod$-R$ fibre lemma and the fact that $S$ surjects on $G/H$ we have that $G/H$ acts nilpotently on the $\mathbb{F} _p-$homology of $BH$ and therefore the map $BG\rightarrow BG\pcom$ is a mod$-p$ equivalence if and only if $BH\rightarrow BH\pcom$ is a mod-$p$ equivalence and therefore $BG$ is $p-$complete if and only if $BH$ is $p-$complete. This implies that $BG$ is $p-$good if and only if $BH$ is $p-$good.
$\Box$
\begin{Proposition}[\cite{BLO3} Proposition 4.4]
The classifying space of a discrete $p-$toral group is $p-$good.
\end{Proposition}
\begin{Remark}
This shows that groups with infinite Sylow $p-$subgroups can be good and bad for the prime $p$.
\end{Remark}
\section{Completions and localisations at different primes and $\mathbb{Z}_{(J)}$}
We illustrate that spaces can be good and bad at different primes in any arbitrary combination. The special case answers a question concerning the compatibility of completions at different primes raised during the talk of Bob Oliver on $p-$local homotopy theory at the Mathematical Sciences Research Institute Berkeley at the beginning of the algebraic topology program in  the spring 2014.\\[0.3cm] Let $\{p_i\}_{i\in I}$ be an arbitrary 
collection of primes. For every prime $p$ let $G_p$ be one of the groups as defined in Proposition \ref{groups}. Let  $\{H_j\}_{j\in\mathcal{J}}$ be a collection of finite $p-$groups,  $p\notin \{p_i\}_{i\in I}$, and not more than one element for each $p\notin \{p_i\}_{i\in I}$ and including the possibility that $H_j$ is trivial. Let $G$ be any arbitrary combination of direct and free products of the groups $G_{p_i}$ and $H_j$, for example $G:=(\underset{i\in I}{\prod  }G_{p_i})\times(\underset{j\in \mathcal{J}}{\prod}H_{j})$ or $G:=(\underset{i\in I}{*}G_{p_i})*(\underset{j\in \mathcal{J}}{*}H_{j})$. Assume that the group $G$ has a finite index normal subgroup generated by elements of $q'-$order for all primes $q\in \mathcal{I}$.
\begin{Theorem}
The classifying space $BG$ 
 is $p_i-$bad for $i\in I$ and $p-$good for all other primes $p$.
\end{Theorem}
\underline{Proof:} In the first case we have $G:=(\underset{i\in I}{\prod  }G_{p_i})\times(\underset{j\in \mathcal{J}}{\prod}H_{j})$. Then for the classifying space we have $BG\simeq
B(
(\underset{i\in I}{\prod}G_{p_i})\times(\underset{j\in \mathcal{J}}{\prod}H_{j}))
\simeq\underset{i\in I}{\prod} BG_{p_i}\simeq BG_{p_i} \times (\underset{j\in I,i\neq j}{\prod} BG_{p_j}) \times ( \underset{j\in\mathcal{J}}{\prod} B H_j)$ for all $i\in I$ and $BG_{p_i}$ is $p_i-$ bad for all $i\in I$ [Proposition \ref{pgoodpbad}]. Lemma \ref{pbad} implies that $BG$ is $p_i$-bad for all $i\in I$. For all other primes $p\notin \{p_i\}_{i\in I}$ we have $G=(\underset{i\in I}{\prod}G_{p_i})\times(\underset{j\in \mathcal{J}}{\prod}H_{j})$ has a finite index normal subgroup generated by all elements of $p'-$order and therefore \cite[ Proposition 1.11]{AKO} the classifying space $BG$ is $p-$good. If we have $G:=(\underset{i\in I}{*}G_{p_i})*(\underset{j\in \mathcal{J}}{*}H_{j})$ we have $BG\simeq
B(
(\underset{i\in I}{*}G_{p_i})*(\underset{j\in \mathcal{J}}{*}H_{j}))
\simeq  (\underset{j\in I}{\bigvee} BG_{p_j}) \vee ( \underset{j\in\mathcal{J}}{\bigvee} B H_j)\simeq BG_{p_i} \vee (\underset{j\in I,i\neq j}{\bigvee} BG_{p_j}) \vee ( \underset{j\in\mathcal{J}}{\bigvee} B H_j)$. This space is $p_i-$bad for all $i\in I$ and  $p-$good for all other primes $p$. The proof of the result for the classifying space of an arbitrary combination of direct and amalgamated products of the groups $G_{p_i}$, $H_j$ is totally analogous.  $\Box$ 
\begin{Remark}
For spaces with finite homology or homotopy groups the completion with respect to the integers is the product of the respective $p$-completions.
\end{Remark}
\begin{Remark}
The question whether the space $BG$ constructed above is good or bad for the integers in the arbitrary case remains open. The question whether a space can be good and bad for arbitrary solid rings is solved negatively in a seperate manuscript \cite{RgoodRbad}. The results of the following section show that the methods of this article do not extend canonically to the case of arbitrary solid rings.
\end{Remark}
\begin{Remark}
For nilpotent spaces the $p-$completion determines the $\mathbb{Z}[p^{-1}] = \mathbb{Z}_{(p)}-$localisation \cite{BK}.
\end{Remark}
\section{Restrictions of
mod-$R$ Fibre Lemma to $\mathbb{Z}_{(J)}$-completions}
The following shows the limitations of the mod $R-$fibre lemma  \cite[Lemma 5.1]{BK} for proving that a space is bad in the sense of Bousfield-Kan concerning completions with respect to the integers. We believe this result is known to specialists even though there seems to be no proof in the literature.
\begin{Theorem}
There is no nontrivial nilpotent action of a finite group on a free abelian group.
\end{Theorem}
\underline{Proof:} We use induction on $n$ to show that if $A=\mathbb{Z}^n$ is a nilpotent $G-$module then it is trivial. If $m=0$ there is nothing to prove, so assume the result holds for all $m<n$. The submodule $A^G=H^0(G,A)$ is not trivial so $A/A^G\cong \mathbb{Z}^m\oplus F$ where $m<n$ and $F$ is finite. If $x\in A$ is a preimage of some $\overline{x}\in F$ then $kx\in A^G$ and since $A$ is torsion-free, $kx=g(kx)=kg(x)$ implies that $g
(x)=x.$ This shows that $F=0$. If $m=0$ then $A$ is a trivial $G-$ module as needed. So assume $m\geq 1$. Since $m<n$, then induction hypothesis implies that $A/A^G$ is a trivial $G-$module, hence $G$ acts via a homomorphism into $\{ \begin{pmatrix}I&*\\ 0 &I \end{pmatrix}\}\cong  \mathbb{Z}^{m(n-m)}$, and therefore the action must be trivial. $\Box$
\begin{Remark}
 The proof extends to $\underset{i\in I}{\prod }{\mathbb{Z}[J^{-1}_i]}$ for arbitrary sets of primes $J_i$ as soon as $|G|$ is coprime to all elements of $J_i$ for all $i\in I$ because the representation theory in this case is semi-simple.
\end{Remark}
\begin{Remark}
 In the case of the integral homology of the fibre $F$ of the fibration associated with the group extension $F\rightarrow C_p*C_p\rightarrow C_p\times C_p$ the modules is not rationally irreducible. It has $p-1$ rational summands, each of dim $p-1$.
\end{Remark}
\begin{Remark}
The same module over the complexes consists of
the sum of all the one dimensional characters that are nontrivial on each of the
two copies of $C_p$.
\end{Remark}
\begin{Remark}
The
 index set $I$ 
is arbitrary and the result extends to abelian or nilpotent groups because of the classification theorem of abelian groups.
\end{Remark}

\section{Completions, localisations and $R-$homological structures
}
Every space is pro-isomorphic to a tower of nilpotent spaces \cite{Farjoun2} ruining the hope that good or bad can be detected from any homological structure. We illustrate this fact in the following discussion. In particular every bad space has the homology of a good space for any solid ring $R$ including the Steenrod algebra structure and the Bockstein spetral sequence (in the cases of $R$ where they exist).\\[0.3cm]
It
 is impossible to detect from the $\mathbb{F}_p-$homology module structure
whether a
space is good or bad.
\begin{Proposition}
Let $X_1$ be the space $S^1\vee S^2$ and $X_2$ be the space $\mathbb{R}P^2$.
 Then $X_1$ and $X_2$ have isomorphic homology as $\mathbb{F}_p-$modules but
 $X_1$
 is bad for all primes and the integers
and $X_2$ is good for all primes and bad for the integers.
\end{Proposition}
\underline{Proof:} 
The $\mathbb{F}_p-$module structure
is determined by $H_i(X_j;\mathbb{F}_p)\cong \mathbb{F}_p$ for $i\in\{0,1,2\}$ and $j\in\{1,2\}$ and $0$ otherwise. The coalgebra structure is different however. The space $X_1$ is bad
for primes and the integers \cite[Theorem 10.1]{Bousfield}
and the space $X_2$ is good for primes and bad for the integers \cite{BK}.
$\Box$\\[0.3cm]
It follows from Lemma \ref{pbad} that the wedge of any space with a bouquet of spheres is bad for all primes as soon as the bouquet consists of at least two spheres and one of them is $S^1$. The Steenrod algebra acts trivially on spherical classes. This motivated the conjecture that if for a space $X$ every element in the reduced homology is either in the source or the target of a nontrivial dual Steenrod operation it might be good. This however is false.
It is not possible to tell from the Steenrod algebra structure whether a space is $p-$good or $p-$bad in the sense that even if every element in the reduced homology module $\overline{H}_*(X;\mathbb{F}_2)$ is
either in the source or in the target of a nontrivial dual Steenrod operation the space can be $2-$good as in the case of $\mathbb{R}P^{\infty}\vee\mathbb{R}P^{\infty}\simeq B(C_2*C_2)$ or $2$-bad as for $\mathbb{R}P^{\infty}\vee\mathbb{R}P^{\infty}\vee\mathbb{R}P^{\infty}\simeq B(C_2*C_2*C_2)$. In fact we have that it is impossible to tell from any $R$-homological structure at all whether a space is $R$-good or $R$-bad for any solid ring $R$.\\[0.3cm]


Dr.\ Nora Seeliger PhD, Department of Mathematics and Statistics, Room B9, Bailrigg Campus, Lancaster University, Lancaster, LA1 4YF, Lancashire, United Kingdom. Email: nora@lancaster.ac.uk.

\begin{thebibliography}{breitestes Label}
\bibitem{AKO} M.\ Aschbacher, R.\ Kessar, B.\ Oliver, \textit{Fusion systems in algebra and topology}, London Mathematical Society Lecture Note Series: 31, Cambridge University Press, 2011.
\bibitem{Bousfield} A.\ K.\ Bousfield, \textit{On the $p-$adic completions of nonnilpotent spaces}, Transactions of the American Mathematical Society, Volume 331, Number 1, May 1992, 335--359.
\bibitem{BK} A.\ K.\ Bousfield, D.\ M.\ Kan, \textit{Homotopy Limits, Completions and Localizations}, Springer Lecture Notes in Mathematics, Springer-Verlag, Berlin, Heidelberg, New York, 1972.
\bibitem{BK2} A.\ K.\ Bousfield, D.\ M.\ Kan, \textit{The core of a ring}, J.\ Pure Applied Algebra 2 (1972), 73--81.
\bibitem{BLO3} C.\ Broto, R.\ Levi, B.\ Oliver, \textit{Discrete models for the $p$-local homotopy theory of compact Lie groups and $p$-compact groups}, Geometry and Topology 11 (2007) 315--427.
\bibitem{CurtisReiner} C.\ W.\ Curtis, I.\ Reiner, \textit{Representation Theory of Finite Groups and Associative Algebras}, Interscience Publishers, New York, 1962.
\bibitem{Dickson} L.\ E.\ Dickson, \textit{Linear Groups with an
Exposition of the Galois Field Theory}, (1901).
\bibitem{Farjoun} E.\ D.\ Farjoun, \textit{Two Completion Towers for Generalized Homology}, Contemporary Mathematics Volume 265 (2000).
\bibitem{Farjoun2} E.\ D.\ Farjoun, \textit{Pro-nilpotent representation of homology types}, Proceedings of the American Mathematical Society, Volume 38, Number 3, May 1973, 657--660.
\bibitem{Hall} M.\ Hall, \textit{The theory of groups}, Macmillan (1959).
\bibitem{IvanovMikhailov1} S.\ O.\ Ivanov, R.\ Mikhailov, \textit{On a problem of Bousfield for metabelian groups}, Adv.\ Math.\ 290 (2016), 552--589. 
\bibitem{IvanovMikhailov2} S.\ O.\ Ivanov, R.\ Mikhailov, \textit{On lenghts of $H\mathbb{Z}-$localization towers}, preprint.
\bibitem{Ian+Radu} I.\ Leary, R.\ Stancu, \textit{Realising fusion systems}, Algebra Number Theory 1.1, (2007), 17--34.





\bibitem{Malcev} A.\ L.\ Malcev, \textit{Nilpotent groups without torsion}, Jzv. Akad. Nauk. SSSR
, Math. 13 (1949), 201--212.
\bibitem{MartinoPriddy} J.\ Martino, S.\ Priddy, \textit{Unstable homotopy classification of $BG\pcom$}, Math.\ Proc.\ Camb.\ Phil.\ Soc.\ 137 (2004) 321--347.
\bibitem{Massey} William S.\ Massey, \textit{Algebraic topology: an introduction}, Graduate Texts in Mathematics, Springer-Verlag, New York, 1977.
%
\bibitem{Quillen} D.\ G.\
 Quillen, \textit{Rational homotopy theory}, Annals of Math.\ 90
 (1969), 205--295.
\bibitem{Robinson1} G.\ Robinson, \textit{Amalgams, blocks, weights, fusion systems and finite simple groups}, Journal of Algebra 314 (2007), 912--923.
\bibitem{gmffs} N.\ Seeliger, \textit{Group models for fusion systems}, Topology and its Applications 159 (2012), no. 12, 2845--2853.
\bibitem{RgoodRbad} N.\ Seeliger, \textit{A few examples of $R-$good and $R-$bad classifying spaces}, arXiv:1703.05754.
\bibitem{Sullivan} D.\ Sullivan, \textit{Geometric topology, part I: localization, periodicity and Galois symmetry}, MIT (1970).
\bibitem{Serre} J.\ P.\ Serre, \textit{Trees}, Springer-Verlag Berlin Heidelberg New York, 1980.
\end{thebibliography}
\end{document}